\documentclass[]{interact}

\usepackage{epstopdf}

\usepackage{natbib}
\bibpunct[, ]{(}{)}{,}{a}{}{,}%
\def\bibfont{\small}%
\def\bibsep{\smallskipamount}%
\def\bibhang{24pt}%
\defcitealias{CDCcost}{CDC 2018a}
\defcitealias{CDCepidemiology}{CDC 2018b}
\defcitealias{WHOincubation}{WHO 2018}
\defcitealias{WHOorigin}{WHO 2015}
\defcitealias{world2014ebola}{WHO 2014}

\theoremstyle{plain}

\theoremstyle{definition}

\theoremstyle{remark}

\usepackage[inline]{enumitem}
\usepackage{amsmath}
\DeclareMathAlphabet\mathbfcal{OMS}{cmsy}{b}{n} 
\usepackage{amssymb}
\usepackage{kpfonts} 
\usepackage{graphicx}
\usepackage{adjustbox} 
\usepackage{longtable}

\usepackage{nicefrac}
\usepackage{multirow}
\usepackage{placeins} 
\usepackage{chngcntr}
\usepackage{float}
\usepackage{subcaption}            
\usepackage[onehalfspacing]{setspace}

\begin{document}

\title{Data-Driven Infectious Disease Control with Uncertain Resources}

\author{
\name{Ceyda Yaba Best\textsuperscript{a}, Amin Khademi\textsuperscript{a}, and Burak Eksioglu\textsuperscript{b}\thanks{CONTACT B. Eksioglu.
Email: burak@uark.edu}}
\affil{\textsuperscript{a}Department of Industrial Engineering, Clemson University, Clemson, SC \\ \textsuperscript{b}Department of Industrial Engineering, University of Arkansas, Fayetteville, AR} }

\maketitle
\begin{abstract}
We study a resource allocation problem for containing an infectious disease in a metapopulation subject to resource uncertainty. We propose a two-stage model where the policy maker seeks to allocate resources in  both stages where the second stage resource is random. Instead of a system of nonlinear differential equations that governs the epidemic trajectories in the constraints of the optimization model, we use a data-driven functional form to model the cumulative number of infected individuals. This flexible data-driven modeling choice allows us to transform the optimization problem to a tractable mixed integer linear program. Our flexible approach can handle an online decision making process, where the decision makers update their decisions for opening treatment units and allocating beds utilizing the new information about the epidemic progress. We utilize a detailed simulation model, validated by real data from the 2014 Ebola epidemic in Sierra Leone. Our results show that our policies produce about 400 fewer number of infected individuals in Sierra Leone compared to the policies applied during the actual epidemic. We also provide a detailed comparison of allocation policies generated by our optimization framework which sheds light on the optimal resource allocation in different regions.
\end{abstract}


\vspace{2in}

\section{Introduction}\label{intro}
Emerging and re-emerging infectious diseases are continuously posing a threat to public health worldwide. It is estimated that more than $25\%$ of deaths worldwide are attributed to an infectious disease \citep{Fauci2005}. Specifically, in low income countries the top three causes of death are related to infectious diseases \citep{WHOdeath}. With the increasing rate of urbanization in low income countries, challenges for designing and operating reliable healthcare infrastructures have become a major concern for governments \citep{Neiderud2015}. For example, in Sierra Leone $\$224$ per capita is spent on health compared to $\$9403$ in the United States \citep{whoSLE, whoUS}. Also, \cite{Raven2018} reported that there are only two doctors available per 100,000 people in Sierra Leone. In addition, health systems have been extremely fragile due to internal and international conflicts to cope with large scale epidemics such as the 2014 Ebola outbreak.

Our study is motivated by unique operational challenges faced by decision makers during large scale epidemics such as the 2014 Ebola epidemic. After the initial lack of proper response combined with poor healthcare infrastructure, international organizations formed a more coordinated relief effort to control the Ebola epidemic in Sierra Leone and the other affected countries, Guinea and Liberia \citep{Ross2017}. Due to the fact that local hospitals were lacking essential infection control equipment and sterile isolation environment, coupled with fear of Ebola among healthcare workers \citep{Gostin2014}, some health organizations and local nonprofits started building Ebola Treatment Units (ETUs) in Sierra Leone \citep{Shin2018}. Opening and maintaining of ETUs provide essential care for the patients, while isolating infected individuals from the rest of the community and imposing strict guidelines \citep{WHOtreatment}.  However, the challenge remains in allocating limited budgets for isolation efforts.

The funding need was estimated by the World Health Organization (WHO), and funds were received from different donors including non-governmental organizations and direct government support \citep{UNresources}. However, due to the nature of these donations there is uncertainty involved in the amount of funds that the decision makers will finally receive. Therefore, a natural question is how to plan for controlling an infectious disease when the amount of funds that is expected through donations is subject to uncertainty. We develop a flexible framework to address this question. In particular, we adopt a two-stage optimization framework in which the known and random budget is used for opening and maintaining ETUs in the first and second stage, respectively. Our framework allows for other types of control mechanisms in the formulation through various parameters. However, we focus on isolation strategies via ETUs. 

In order to account for disease dynamics, we use a novel approach that takes advantage of the shape of the curves that correspond to the cumulative number of infected individuals. If one follows a standard approach, the infection is governed by a system of nonlinear differential equations. Including the discretized versions of these differential equations will make our two-stage formulation nonlinear and intractable. Therefore, we use the observation that cumulative number of infected individuals follows an \textit{s}-shaped curve as a function of time \citep{Brandeau2005} to propose a data-driven mechanism. In particular, this observation allows us to replace the system of nonlinear differential equations for prediction of the cumulative number of infected individuals with a function whose parameters are estimated by a validated simulation model of the system of nonlinear differential equations. Based on the observations in the literature and our simulation, we fit logistic curves to obtain a closed form function of time which depends on the initial state and corresponding bed allocations. Moreover, our extensive simulation experiments show that fitting a linear regression for one of the parameters of the functional form produces high quality results. Therefore, using such a data-driven approach, the resulting two-stage optimization problem becomes a two-stage stochastic program (TSSP) with linear constraints and integer variables.
 
Clearly, the quality of produced solutions depends on the accuracy of parameter estimation. To estimate the parameters of our logistic functions, we develop a detailed simulation model using data from the 2014 Ebola epidemic in Sierra Leone. We use an online approach to mimic the real decision making process. We first assume that there is a limited amount of data available in the first stage (earlier stages of the epidemic). The decision maker uses the validated simulation model with respect to the data available on hand to generate trajectories of the epidemic, and estimates the parameters of the logistic curves using these trajectories. Then, she solves an optimization problem to find actions for opening treatment units and bed allocations. In the beginning of the second stage, after more data is observed, the decision maker updates the information by re-calibrating/re-validating her simulation model to find new parameters for the logistic curves. This procedure allows the decision maker to update the second stage decisions based on new information. When we adopt this approach to estimate the parameters in our TSSP model for opening treatment units and bed allocations, our results show that the cumulative number of infected individuals is improved by about $3\%$ compared to the real bed allocations used in Sierra Leone. This improvement translates to averting around 400 infections. We also produce solutions retrospectively by using all the historical data which is a common practice in data-driven approaches. By using the whole data, we find the best decisions that could have been taken if the decision makers had perfect information on the disease progression. In this setting, our results show an improvement up to $5\%$ is possible when compared to what had taken place during the real response in Sierra Leone. Our analysis suggests that the value of information is significant when making decisions for a large scale epidemic control. 

This paper provides a detailed analysis for a decision maker when she does not have access to vast amounts of data when facing a large scale epidemic with limited resources. The data-driven optimization literature takes advantage of having vast amounts of historical data (see \cite{Bertsimas2013}). However, our online decision making process allows us to exploit validated simulation models when there is limited data available, and update the information after the system is observed for some more time. By using this approach, the decision maker can take advantage of updated information to improve on her decision making process. 

In summary, we make the following contributions: \begin{enumerate*}[label=\arabic*.]
\item We formulate the problem of resource allocation for containing an epidemic in a metapopulation (a collection of possibly interacting populations) when resources for the second stage are uncertain.
\item We introduce a data-driven approach to solve the optimization problem by which an intractable formulation is transformed to a mixed-integer linear program. This framework is flexible in that the process can be implemented in an online fashion.
\item We develop a validated simulation model to estimate the parameters of our data-driven approach.
\item Our results show that significant improvement is possible by optimizing actions regarding opening treatment units and bed allocations. Also, we provide insights on the optimal resource allocation in different regions of Sierra Leone. 
\end{enumerate*} The rest of this paper is organized as follows: Section \ref{literature} reviews the recent literature. Section \ref{model} describes our modeling process for the resource allocation problem. We provide numerical results in Section \ref{result}. Section \ref{conclusion} concludes with insights and future research directions.

\section{Literature Review}\label{literature}
This work focuses on data-driven resource allocation under uncertainty for infectious disease control in a metapopulation. We briefly explain four main streams of literature related to our study and discuss our contributions.

Resource allocation problems are widely studied in infectious disease control literature. There are many different analysis techniques that attract the attention of researchers from operations research to epidemiology. For example, \cite{Malvankar-Mehta2011} study a two-level resource allocation model where the upper level decision maker allocates resources to lower level decision makers who distribute these resources to local programs. \cite{Ekici2014} model a facility location problem for food distribution in an influenza epidemic where they estimate the need for food by an agent based simulation model of the epidemic. \cite{Juusola2015} suggest a linear model to allocate resources to a mix of HIV prevention and treatment efforts. 
\cite{khademi2015price} measure the price of reallocating HIV treatment in resource-limited settings.
All these studies present results for epidemics in a single population with deterministic resources. However, we consider a metapopulation with interacting populations and random resources. 

\cite{Mbah2011} use an optimal control problem for two coupled populations to study efficiency and equity. \cite{Kasaie2013} consider a resource allocation problem for multiple interacting populations using a simulation/optimization framework. \cite{Alistar2013} develop a model to allocate HIV prevention and treatment resources while presenting results for a single population, multiple independent populations as well as two interacting populations. A recent study by \cite{Long2018} also considers a metapopulation with interacting populations and provide spatial bed allocations subject to limited budget for the 2014 Ebola epidemic. These models are useful to a decision maker facing an epidemic control problem, however they assume that the resources to be allocated are known ahead of time. In our work, we consider uncertainty in the budget for resource allocation, which is different from the aforementioned models.

Another stream of literature that is related to our work is decision making under budget or supply uncertainties. \cite{Natarajan2014} formulate an inventory control problem for humanitarian supply chains with uncertain budget installment periods. They assume that the total amount of the budget is known, however the decision makers do not know the number of installments and the amount in each installment. There are also many studies that focus on the uncertainty in the supply rather than budget. For example, \cite{Oezaltin2011} formulate a multi-stage stochastic program for influenza vaccination design where the supply is uncertain due to several factors in the production process.
A major difference between our setup with previous models studied in the literature lies in the fact that our objective function is based on the number of infected individuals in the final period, which depends on the random budget affecting a system of nonlinear differential equations. Our contribution is to develop a data-driven integrated simulation-optimization framework to make the problem tractable. Our data-driven estimation of epidemic progression is also a novel approach, which is incorporated in an optimization framework to address budget uncertainty.     
Finally, while other models provide useful results in their respective settings for supply uncertainty, our model focuses on the budget uncertainty for a resource allocation model in an epidemic over a metapopulation, which is a novel application. The models that assume uncertainty for infectious disease control in resource allocation usually focus on the uncertainty to measure the effectiveness of control strategies, e.g., see \cite{Tanner2008} and \cite{Yarmand2014}.

Recently, data-driven optimization has gained attention in operations research and management science communities. \cite{Misic2019} provide  a detailed survey on recent data-driven optimization models. Here, we focus on the literature for data-driven models in healthcare settings. \cite{Bertsimas2013} use historical data to find a scoring based priority policy for kidney transplantation in the U.S. They use a matching problem where the parameters of the optimization model are found by historical data. They then fit a regression model to find weights on the characteristics of the patients and organs that are important to decision makers in a ranking sytem for efficiency and fairness. \cite{Bertsimas2016} utilize regularized regression models to estimate the outcomes of clinical trials and formulate an optimization problem for selecting drugs using estimates from historical data. \cite{Rath2018} formulate a TSSP for staff planning in hospitals where demand for anesthesiologists is uncertain, and the cost parameters are estimated using historical data. The models mentioned here assume a vast availability of data from the past and estimate the parameters using this information retrospectively. However, in our approach, we adopt an online decision making process where a decision maker updates the information on the disease progression with respect to observed data. We also provide an analysis under perfect information which is similar to the aforementioned studies. 

The last stream of research focuses on simulation models designed for the 2014 Ebola epidemic. Most of these studies use deterministic compartmental models for their analysis such as \cite{Althaus2014}, \cite{Fisman2014}, and \cite{Lewnard2014}. \cite{Rivers2014}, as one exception, use a stochastic compartmental model for the national data of Liberia and Sierra Leone. \cite{Pandey2014} use a continuous time stochastic model for the data of the Ebola epidemic in West Africa. These models only focus on national level data, but we consider a metapopulation based on geographic locations and interacting populations, coupled with change of individual behavior. We use a stochastic compartmental model for validation based on data availability and use this information to aid the decision making process.

\section{Problem Formulation}\label{model}
As discussed in Section \ref{intro}, we study the problem of containing an infectious disease in a metapopulation when resources are limited and uncertain. The motivation for the work is based on challenges faced by decision makers in the 2014 Ebola outbreak in West Africa. In particular, we consider two special features: \begin{enumerate*}[label=(\roman*)]
\item metapopulation: the epidemic started in a remote forest region in the Gueckedou district in Guinea and spread to Liberia and Sierra Leone \citep{WHOorigin}; therefore our setting includes the interaction among populations such as cities;
\item limited and uncertain resources: the epidemic happened in areas that are limited in resources and more importantly they were uncertain, a novel feature that we study in this work. More specifically, in the beginning of the epidemic, the decision makers have to open treatment units and then assign beds to those treatment units. After the initial phase, the epidemic had become a global issue and additional resources became available from the relief efforts and donations from non-profit organizations and governments. However, these extra resources are subject to uncertainty, especially in the early phases of the epidemic when the decision makers seek to open treatment units in different geographical locations.
\end{enumerate*} 

Motivated by such challenges, we study the containment of an infectious disease where the decision maker has a limited budget in the first stage to open treatment units in possibly different geographical locations and allocate beds to these units. In the second stage, when the uncertain extra resources are revealed, the decision maker has to allocate beds to treatment units already opened in the first stage. Note that, for the second stage, we only allocate additional beds to existing treatment units, because in the motivating application, setting up a new treatment unit is challenging and time consuming. However, the mathematical formulation that we propose is flexible to consider opening treatment units in the second stage as well. In both stages, we have a system of nonlinear differential equations that describe how the number of individuals in every compartment evolve over time (see Section \ref{simulation} for the details of the compartmental model). Therefore, the problem formulation becomes intractable. To overcome this, we propose a novel approach to make the problem a tractable TSSP with linear constraints and integer variables. The idea behind our proposal is to represent the cumulative number of infected individuals by a logistic curve, which is widely used for studying infectious disease dynamics \citep{brauer2012}. In particular, we let $I_n(t, \mathbf{I}(0), \mathbf{m}) = \frac{K_n(\mathbf{I}(0), \mathbf{m})}{1+e^{-(a_n t + b_n)}}$ denote the expected cumulative number of infected individuals in population $n$ at time $t$, where $\mathbf{I}(0)$ denotes the vector of initial number of infected individuals in each population, and $\mathbf{m}$ denotes the vector of resources assigned to every population. Moreover, we assume that the function $K_n(\mathbf{I}(0), \mathbf{m})$ is linear in its arguments, which makes the formulation a linear TSSP (see Section \ref{s_shaped} for details of the assumptions and justifications). A natural question is  ``how can we estimate the parameters involved in $I_n(t, \mathbf{I}(0), \mathbf{m})$?" If the decision maker uses the system of nonlinear differential equations, the parameters to use in the model would be death rates, infectivity rates, mixing patterns, and so on. In that case, the decision maker has to somehow estimate these parameters by data and try to solve an extremely difficult optimization problem. Instead, we propose using $I_n(t, \mathbf{I}(0), \mathbf{m})$ and estimate the parameters using a detailed large-scale simulation model of the epidemic and solve a tractable optimization problem. In particular, we propose an online procedure to mimic the decision making process in real life. That is, the decision maker creates a simulation model of the epidemic, calibrates and validates it by the limited data that is available in the beginning of the first stage. This simulation is used to generate trajectories of the epidemic for the future, which are then used to estimate the parameters of $I_n(t, \mathbf{I}(0), \mathbf{m})$. After observing the state of the system, the decision maker has an option to update the future decisions when she observes the second stage uncertainty. 
This is done by re-calibrating/re-validating the simulation model to generate future epidemic trajectories in the beginning of the second stage and solve the optimization problem for the second stage.
Next, we describe the details of the three models: optimization model, estimation procedure, and simulation development.

\subsection{The Proposed Model}\label{stochastic_program}
We define a discrete and finite set $\mathbf{\Omega}$ for all possible amounts of resources that might become available in the second stage and its probability mass function by $p(\omega)$ for $\omega \in \mathbf{\Omega}$. We divide the planning horizon into two; the first stage up to some (known) time $\tau$ where the initial budget is used to open ETUs and allocate beds, and the second stage where the additional budget is revealed and used to allocate additional beds. Recall that, we do not consider opening extra treatment units in the second stage because the construction time is long and the epidemic must be contained in a short time period. 
For example, \citet{janke2017beyond} reported that the construction time for the installation of an ETU in Monrovia was about 4 months.
Nonetheless, this assumption can be relaxed. The sets, parameters, and the decision variables for the TSSP are given below.

\noindent \textbf{Sets}\\	
$\mathbf{\Omega}$: Set of all possible outcomes for the random additional budget,\\
$\mathcal{N}$: Set of populations.\\

\noindent \textbf{Parameters}	\\
$I_n(0)$: Initial number of infected individuals in population $n \in \mathcal{N}$ ($\mathbf{
I}(0)$ is the vector of initial infected individuals for the metapopulation),\\
$B$: Available budget in the first stage,\\
$o_n$: Cost of opening an ETU in population $n \in \mathcal{N}$,\\
$h_n$: Cost of maintaining a bed in an ETU in population $n \in \mathcal{N}$,\\
$r_n$: Fairness rate of population $n \in \mathcal{N}$,\\
$\tau$: Time period that the first stage ends,\\
$\tilde{B}(\omega)$: Additional budget for the second stage under scenario $\omega$,\\
$p({\omega})$: Probability of realizing scenario $\omega$,\\
$T$: Time period that the second stage ends.\\

\noindent \textbf{Decision Variables:}\\
$y_n$: $1$, if a treatment unit is opened in population $n$; $0$, otherwise, \\
$m_n$: Number of beds allocated to population $n$ in the first stage ($\mathbf{m}$ is the vector of beds allocated in the metapopulation), \\
$I_n(\tau, \mathbf{m}, \mathbf{I}(0))$: Cumulative number of infected individuals at time $\tau$ in population $n$ ($\mathbf{I}(\tau)$ is the corresponding vector for the metapopulation),\\
$\tilde{m}_n(\omega)$: Number of beds allocated to population $n$ in the second stage under scenario $\omega$ ($\mathbf{\tilde{m}}(\omega)$ is the corresponding vector for the metapopulation $\omega$),\\
$I_n(T, \mathbf{m}, \mathbf{I}(\tau), \mathbf{\tilde{m}}(\omega))$: Cumulative number of infected individuals at time $T$ in population $n$, with respect to scenario $\omega$.

The problem can now be formulated as
\begin{subequations}\label{sp_formulation}
\begin{flalign}
& \min \sum_{n \in \mathcal{N}} \sum_{\omega \in \mathbf{\Omega}} p(\omega) I_n(T, \mathbf{m}, \mathbf{I}(\tau), \mathbf{\tilde{m}}(\omega)) \label{obj}&\\
& \text{s.t.}\notag&\\
& \qquad \sum_{n \in \mathcal{N}} o_n y_n + h_n m_n \leq B, \label{budget1}&\\
&\qquad m_n \leq M y_n, \qquad \forall n \in \mathcal{N}, \label{open1}&\\
&\qquad m_n \leq r_n \sum_{n' \in \mathcal{N}} m_{n'}, \qquad \forall n \in \mathcal{N}, \label{fairness}&\\
&\qquad \sum_{n \in \mathcal{N}} h_n \tilde{m}_n(\omega) \leq \tilde{B}(\omega), \quad \forall \omega \in \mathbf{\Omega}, \label{budget2}&\\
&\qquad \tilde{m}_n(\omega) \leq M y_n, \quad \forall n \in \mathcal{N},   \omega \in \mathbf{\Omega}, \label{open2}&\\
&\qquad m_n+\tilde{m}_n(\omega) \leq r_n \sum_{n' \in \mathcal{N}} m_{n'} + \tilde{m}_{n'}(\omega),\quad \forall n \in \mathcal{N}, \omega \in \mathbf{\Omega}, \label{fairness2}&\\
&\qquad y_n \in \{0, 1\}, \quad \forall n \in \mathcal{N}, \label{binary}&\\
&\qquad m_n, \tilde{m}_n(\omega) \in \mathbb{Z}_{+}, \quad \forall n \in \mathcal{N}, \omega \in \mathbf{\Omega}, \label{integer}&   
\end{flalign}
\end{subequations} 

\noindent where the objective is to minimize the sum of expected cumulative number of infected individuals at time $T$ with respect to different scenarios; $I_n(T, \mathbf{m}, \mathbf{I}(\tau), \mathbf{\tilde{m}}(\omega))$ is a function of time, first and second stage bed allocations, and the cumulative number of infected individuals at time $\tau$ in the metapopulation; constraint \eqref{budget1} is the budget constraint for opening treatment units and allocating/maintaining beds in the first stage. Constraints \eqref{open1} ensure that no bed is allocated if a treatment unit is not opened in a population. We also include a fairness constraint in \eqref{fairness}, such that we impose an upper bound on the number of allocated beds to a population. Constraint \eqref{budget2} is the budget constraint for additional bed allocations in the second stage with respect to random extra budget. Constraints \eqref{open2} are similar to \eqref{open1}, and ensure that no additional beds are allocated if there is no treatment unit. We also introduce the fairness constraints for the second stage bed allocations in \eqref{fairness2}. 

Constraints \eqref{fairness} provide a flexible framework to incorporate fairness in the model. 
In essence, they apply an upper bound ($r_n$) on the percentage of beds allocated to each population, which is easy to understand and implement. This framework enables policy makers to set upper bounds $r_n$ based on a variety of factors. For example, one may set $r_n$ to be the percentage of population $n$ (or infected population or susceptible population) to the total population (or infected population or susceptible population) in the metapopulation. This fairness framework borrows ideas from applications in distributing healthcare resources, where the policy maker applies bounds on the percentage of resources allocated to subgroups in the population. For example, in transplanting valuable organs to a heterogeneous patient population on a waiting list, a common fairness framework enforces a lower bound on the percentage of high quality organs allocated to specific patient groups \citep{Bertsimas2013}. Fairness constraints \eqref{fairness2} follow the same lines of reasoning for the entire planning horizon.       

As discussed before, a major challenge is to find a closed form formula for $I_n(T, \mathbf{m}, \mathbf{I}(\tau), \mathbf{\tilde{m}}(\omega))$ and estimate its parameters. In the next section, we describe this procedure and the quality of solutions produced by it.

\subsection{Estimating the Cumulative Number of Infected Individuals}\label{s_shaped}
Logistic curves are often associated with population growth models. In epidemic modeling, they are also used as a tool for studying dynamics of the infectious diseases \citep{brauer2012}. 
In addition to these reasons, we also observe logistic curves from real data and our validated simulation model (see Section \ref{validation}).  Therefore, we use an \textit{s}-shaped logistic curve in the form $\frac{K}{1+e^{-(at+b)}}$. Next, we explain the process for estimating these curves and finding a function which depends on the initial number of infected individuals and the bed allocations, to be used in our TSSP model defined in Formulation \eqref{sp_formulation}. 

We use the time series generated by the validated simulation model to estimate the cumulative number of infected individuals in every population under status quo, where we use the real initial state of the metapopulation and the bed capacities. 
An appealing feature of our model is that because we use simulation we can consider a complex and realistic model for the progression of the epidemic while the optimization module is still tractable. That is, the simulation model is stochastic and incorporates interactions in a system of nonlinear differential equations.

Since we divide the planning horizon into two stages, we estimate two logistic curves for every population for each stage using the simulated time series; one for time period from the beginning until $\tau$, and another from $\tau$ up to $T$. Note that one may use this procedure in an online fashion, i.e., upon availability of more data the decision maker can re-calibrate/re-validate the simulation model, generate new time series, and have an updated version for estimating $I_n(\cdot)$. In particular, we assume that the functional form for the logistic curve for the cumulative number of infected individuals at time $t$ for population $n$ is given by:
\begin{equation}\label{logistic_curves}
I_n(t) = \begin{cases}
\frac{K_n^1}{1+e^{-(a_n^1 t+b_n^1)}}, \text{ for } 0\leq t \leq \tau, \\
\frac{K_n^2}{1+e^{-(a_n^2 (t-\tau)+b_n^2)}}, \text{ for } \tau < t \leq T. 
		 \end{cases}
\end{equation}  

We estimate the parameters in two steps. In the first step, we estimate all the parameters by fitting the curves retrospectively to the real data (or by generating time series using our validated simulation model if data is not available). At this stage, all of $K_n^1$, $a_n^1$, $b_n^1$, $K_n^2$, $a_n^2$, and $b_n^2$ are estimated. Then, we assume that $a_n^1$, $b_n^1$, $a_n^2$, and $b_n^2$ are fixed during an epidemic and do not change by different bed allocations. This assumption is based on observations from running a variety of settings in our simulation model, which show that the change in these parameters do not lead to significantly different end results. Also, having this assumption enables us to formulate the optimization problem as a mixed integer linear program. Note that for estimating the parameters for the second stage, we use $I(\tau)$ as the initial infected population. 

In the second step, we develop a linear regression model for $K_n^1$ and $K_n^2$ where the independent variables are the initial number of infected individuals and the bed capacities for every population. We use lasso regression since the number of independent variables are fairly large ($2|\mathcal{N}|$) and conventional least square estimates lead to low $R^2$ values. In lasso regression, the idea is to minimize the sum of squared errors and $L_1$ norm penalty for the regression coefficients such that the number of non-zero regression coefficients are minimized. In particular, we assume that 
\begin{align*}
K_n^1 &= \rho_{n,0}^1 + \sum_{n' \in N} \rho_{n', 1}^1 I_{n'}(0) + \rho_{n', 2}^1 m_{n'},\\
K_n^2 &= \rho_{n,0}^2 + \sum_{n' \in N} \rho_{n', 1}^2 I_{n'}(\tau) + \rho_{n', 2}^2 (m_{n'} + \tilde{m}_{n'}),
\end{align*}  

\noindent where $\rho_{n,0}^1$, $\rho_{n,1}^1$, $\rho_{n,2}^1$, $\rho_{n,0}^2$, $\rho_{n,1}^2$, and $\rho_{n,2}^2$ are the regression parameters estimated via lasso for population $n$. To that end, we run several simulations by changing the initial populations ($\mathbf{I}(0)$) and bed capacities ($\mathbf{m}$, $\mathbf{\tilde{m}}$). Then, the simulation provides $\mathbf{I}(T)$, which enables us to run lasso regression. Our numerical results show that the regression models have high $R^2$ values and coefficients are significant.
We can now rewrite the decision variables $I_n(\tau, \mathbf{m}, \mathbf{I}(0))$ and $I_n(T, \mathbf{m}, \mathbf{I}(\tau), \mathbf{\tilde{m}}(\omega))$ as follows:
\begin{align*}
I_n(\tau, \mathbf{m}, \mathbf{I}(0)) &= \frac{\rho_{n,0}^1 + \sum_{n' \in \mathcal{N}} \rho_{n', 1}^1 I_{n'}(0) + \rho_{n', 2}^1 m_{n'}}{1+e^{-(a_n^1 \tau + b_n^1)}}, \quad \forall n \in \mathcal{N}, \\
I_n(T, \mathbf{m}, \mathbf{I}(\tau), \mathbf{\tilde{m}}(\omega)) &= \frac{\rho_{n,0}^2 + \sum_{n' \in \mathcal{N}} \rho_{n', 1}^2 I_{n'}(\tau) + \rho_{n', 2}^2 (m_{n'} + \tilde{m}_{n'}(\omega))}{1+e^{-(a_n^2 \tau + b_n^2)}}, \quad \forall n \in \mathcal{N}, \omega \in \mathbf{\Omega}, 
\end{align*} 

\noindent such that our stochastic programming model becomes linear.

\subsection{The Simulation Model}\label{simulation}
The 2014 Ebola epidemic in West Africa is one of the largest Ebola epidemics with 28,000 cases and 11,000 deaths reported in the affected countries \citep{WHO2016a}. The human-to-human transmission was the main reason for the epidemic growth, and there was no cure or vaccination available at the time. Therefore, non-pharmaceutical activities should be put into effect to contain the disease \citep{Pandey2014}. One of these activities  is the establishment of ETUs, where the patients are isolated from the community and provided with care in proper settings \citep{Mallow2018}. With this motivation, we study the dynamics of Ebola by constructing a simulation model to use in our analysis. 

\begin{figure}[h!]
	\centering
	\includegraphics[scale=0.6]{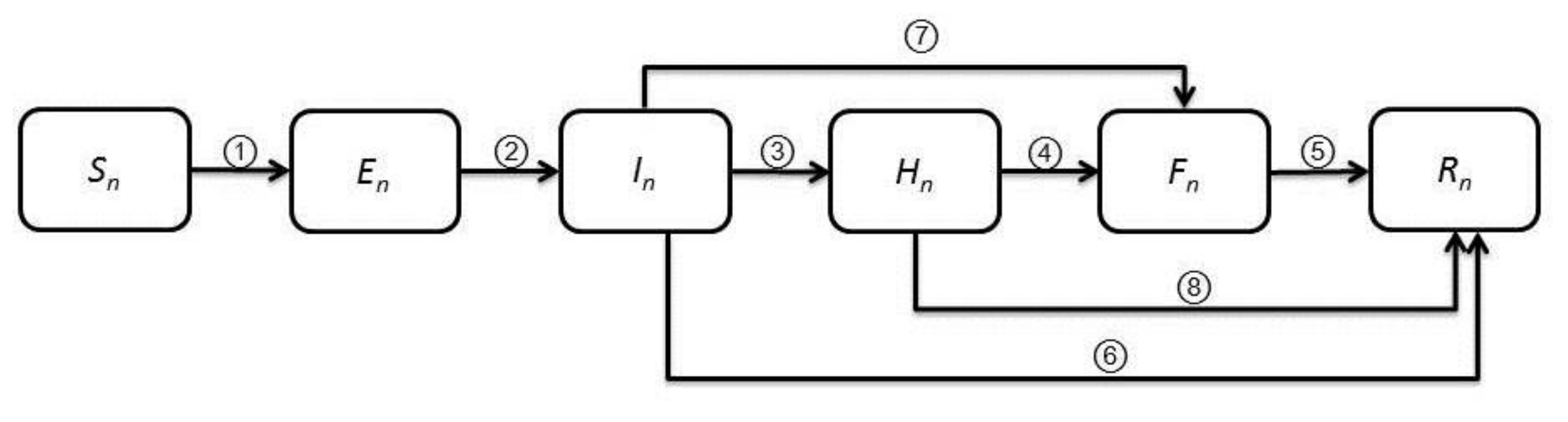}
	\caption{Compartmental model used in the stochastic simulation model for population $n \in \mathcal{N}$.}
	\label{simulation_compartments}
\end{figure}

We develop a stochastic compartmental model to represent the epidemic progression that is adopted from \cite{legrand2006}. 
Our first contribution compared to \citet{legrand2006} is that their simulation is created and calibrated for the whole population, but we consider a metapopulation and model interactions among populations. 
Our second contribution is that we incorporate individual behavior change during the course of an epidemic, which was a significant factor for the 2014 Ebola outbreak.
We calibrate and validate our simulation model accordingly.
There are six compartments in every population; susceptible ($S$), exposed ($E$), infected ($I$), patients in ETUs and hospitals ($H$), patients who die and may continue to infect other susceptible individuals during traditional burial ceremonies ($F$), and patients who do not infect others if they recover or are buried properly ($R$). This compartmental model is appropriate for Ebola, as the incubation period is 2 to 21 days \citep{WHOincubation}, and funerals might be responsible for increase in the reproduction number \citep{Curran2016}. The compartmental model is shown in Figure \ref{simulation_compartments}.

The transition dynamics in population $n \in \mathcal{N}$ are defined as follows: \begin{enumerate*}[label=\arabic*.]
\item A susceptible individual travels to population $n'$ with probability $c_{n, n'}$. During the course of the epidemic, the individuals might change their behavior such that the transmission is reduced. This dampening effect is captured by $e^{-\psi_n t}$ at time $t$, where $\psi_n \ge 0$ is a behavior change coefficient. An individual might become infected upon contact with other infected individuals, patients in treatment units, and the funerals with probabilities $\xi_{n}^I$, $\xi_{n}^H$, and $\xi_{n}^F$, respectively. Thus, a susceptible individual becomes exposed with probability $(e^{-\psi_n t}) \sum _{n' \in \mathcal{N}} c_{n,n'} \frac{\xi_{n'}^I I_{n'} + \xi_{n'}^H H_{n'} + \xi_{n'}^F F_{n'}}{P_ {n'}}$, where $P_{n}$ is the total number of individuals in population $n$.
\item An exposed individual becomes infected with probability $1/\alpha$, where $\alpha$ is the incubation period. 
\item An infected individual is admitted to an ETU with probability $\gamma_{H} \theta$, where $\theta$ is the hospitalization rate and $1/\gamma_H$ is the average duration until admission.
\item A patient in an ETU dies and is buried in a traditional funeral with probability $\gamma_{DH} \delta$, where $\delta$ is the infection related death rate, and $1/\gamma_{DH}$ is the mean duration between admission and death.
\item An infected person who dies is buried and does not infect other susceptible individuals with probability $\gamma_{F}$, with $1/\gamma_F$ as the average time to be buried.
\item An infected individual recovers with probability $\gamma_{I} (1-\theta)(1-\delta)$, where the average time until recovery is $1/\gamma_I$.
\item An infected individual dies and gets buried in a traditional funeral with probability $\gamma_D (1-\theta) \delta$, where the average time until death is $1/\gamma_I$.
\item An individual in an ETU recovers with probability $\gamma_{IH} (1-\delta)$, and the average time until admission and recovery is $1/\gamma_{IH}$. 
\end{enumerate*}
We assume that all the transitions follow binomial distributions, by which the expected values coincide with the deterministic equivalent of the compartmental model.   

We calibrate the simulation model by estimating some of the parameters and using estimates from literature for the others, and each time unit is a day for this study. The definitions and sources of the parameters are given as follows:\\
$P_n$: Total number of individuals in population $n$ \citep{EbolaData},\\
$c_{n,n'}$: Travel rate between populations $n$ and $n'$ \citep{Yang2015},\\
$d_{n, n'}$: Distance between populations $n$ and $n'$ (Google Maps),\\
$\xi_n^I$: Probability of getting infected by contact with infected individuals in population $n$ (calibrated using real data),\\
$\xi_n^H$: Probability of getting infected by contact with infected individuals in ETUs in population $n$ (calibrated using real data),\\
$\xi_n^F$: Probability of getting infected by contact with funerals in population $n$ (calibrated using real data),\\
$\psi_n$: Behavior change coefficient of susceptible individuals for population $n$ (calibrated using real data),\\
$\alpha$: Incubation period (days) \citep{Merler2015},\\
$\theta$: Hospitalization rate \citep{Rivers2014},\\
$\delta$: Death rate \citep{Rivers2014},\\
$1/\gamma_H$: Mean duration from infection to admission to ETUs \citep{Merler2015},\\
$1/\gamma_{DH}$: Mean duration from admission to ETUs to death \citep{Rivers2014},\\
$1/\gamma_{F}$: Mean duration from death to burial \citep{Rivers2014},\\
$1/\gamma_I$: Mean duration from infection to death \citep{Gomes2014},\\
$1/\gamma_{IH}$: Mean duration from admission to ETUs to recovery \citep{Rivers2014}.\\

In calibrating the simulation, we first generate time series for each compartment using an initial set of parameters. We construct a nonlinear optimization problem by using these time series for which the objective function is to minimize the sum of squared errors between the cumulative number of infected individuals generated by the simulation and real data, and the constraints are the discretized versions of the nonlinear differential equations that correspond to the compartmental model. We solve this optimization problem to estimate the parameters of interest as the decision variables, and we continue this process until the difference between two consecutive objective values is smaller than a threshold. We report the calibration and validation errors in Section \ref{validation} where we use the real data from the 2014 Ebola epidemic in Sierra Leone.

\vspace{-4mm}
\section{Numerical Results} \label{result}
In this section, we present numerical contributions by which we shed light on the value of information and optimization in infectious disease control with uncertain resources. The conclusions drawn are based on a validated simulation model of the 2014 Ebola epidemic in Sierra Leone and the results show that improvement regarding patient death is significant. Recall that in order to quantify the value of information in this context, we consider three settings: \begin{enumerate*}[label=\arabic*.]
\item The decision maker has little available data (20 days in our implementation) about the state of the epidemic and has a probability distribution over uncertainty regarding future budget. She creates and validates the simulation and solves the optimization problem, the solution of which provides first and second stage decisions. She then applies these solutions in practice.
\item The decision maker follows the same steps as in Setting 1, except that in the ``middle" of the epidemic when the extra resources become available (60 days in our implementation), she re-optimizes and finds possibly new second stage decisions for implementation. Note that the first stage decisions are irrevocable.
\item The decision maker retrospectively analyzes the epidemic and allocation decisions in the first and second stages. That is, the data is available for the whole period and the goal is to find the best decisions that could have been taken conditioned on the observed epidemic trajectories. Note that this approach is common in data-driven analysis because most analysis is post-implementation, i.e., a dataset becomes available from past behavior of a system and an analysis is run over the data set to draw conclusions. Our approach goes one step further and considers an online process for this application.
\end{enumerate*}

For each of these three settings we report the results of a series of analysis. First, we report the results of simulation calibration/validation. We emphasize that the validation step is important because conclusions are based on the results of the simulation. Second, we describe the procedure to estimate the coefficients of logistic curves for stages one and two and report quantities of interest to validate their quality. Third, using the simulation model, we report the number of infections that could have been averted if our proposed approach had been used. We also report the bed allocation for each setting and conduct sensitivity analysis over the initial budget.  All of the experiments are run on a high performance cluster computer setup in an environment with 2.4 GHz Intel Xeon E5 processor with 20 cores, and 62GB RAM, with Python 3.7.

\subsection{Validation of the Simulation Model}\label{validation}
Since we use the simulation model in three settings as described above, we calibrate and validate it for each setting separately. Data from the 2014 Ebola epidemic in Sierra Leone with 11 districts (two districts are omitted because of the low number of cases in the first 20 days) are used to calibrate the simulation model with  $80\%$ of the available data for training and the rest for validation purposes. Table \ref{calibration_validation} shows the average percent difference between the simulated and real data for calibration and validation periods in each setting in Sierra Leone at the national level. Our validation results show good fit in each setting. Figure \ref{calibration_figures} shows the corresponding time series for cumulative number of infected individuals for each setting. The simulation model reasonably follows the trajectories of the epidemic in all of our settings. Therefore, we can use these simulated curves for our curve fitting process. We report the full results from our validated simulation model for every district in Appendix \ref{districts_validated}.

\begin{table}[H]
	\centering
	\caption{Percent difference between the real and simulated cumulative number of infected individuals.}
	\vspace{0.2cm}
	\large
	{\begin{tabular}{|c|c|c|}
			\hline
			& Calibration & Validation \\ \hline
			Setting 1& $9\%$ & $6\%$ \\ \hline
			Setting 2 & $16\%$ & $8\%$ \\ \hline
			Setting 3 & $11\%$ & $-4\%$ \\ \hline
	\end{tabular}}
	\label{calibration_validation}
\end{table}

\begin{figure}[H]
	\centering
	\begin{subfigure}{\textwidth}
		\centering
		\includegraphics[scale=0.7]{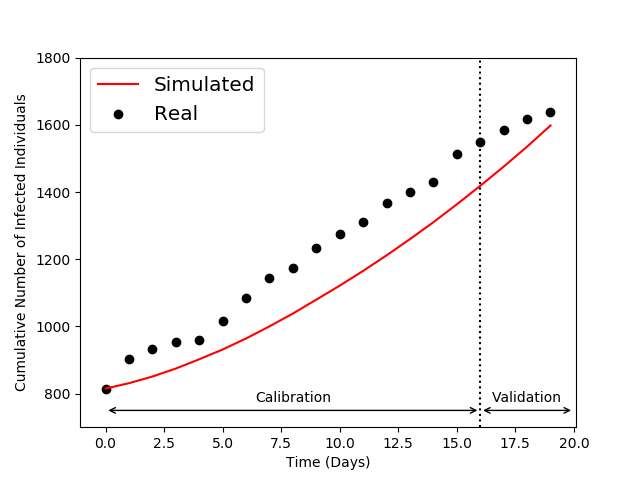}
		\caption{Setting 1}
	\end{subfigure}
\end{figure}

\begin{figure}[H]\ContinuedFloat
	\centering
	\begin{subfigure}{\textwidth}
		\centering
		\includegraphics[scale=0.7]{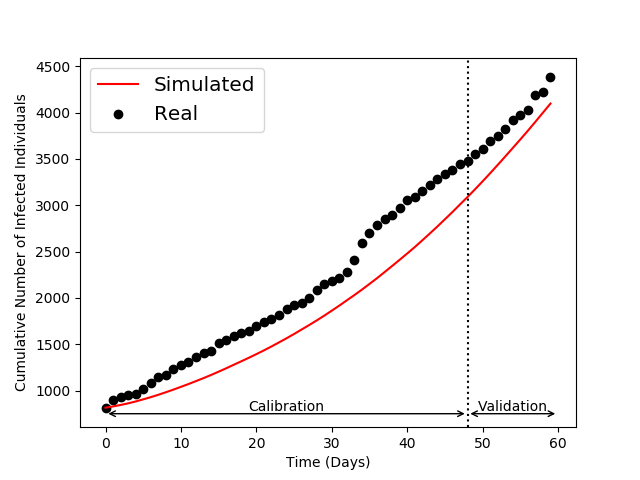}
		\caption{Setting 2}
	\end{subfigure}
	\\
	\begin{subfigure}{\textwidth}
		\centering
		\includegraphics[scale=0.7]{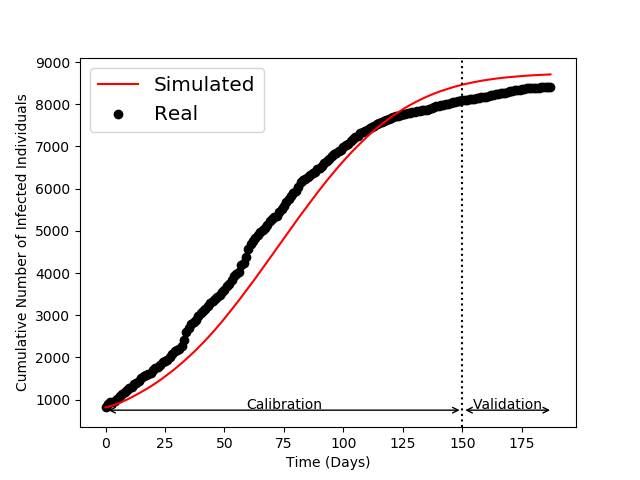}
		\caption{Setting 3}
	\end{subfigure}
	\caption{Calibration and validation of the simulation model with respect to different data availability.} \label{calibration_figures}
\end{figure}

\vspace{-4mm}
\subsection{Functional Estimations for Cumulative Number of Infected Individuals}
In the previous section, we reported the validation results for the simulation model. In this section, though, we report the validation results for fitting the logistic curves in the first and second stages for all three settings. In particular, we report the percent difference between $I(\tau)$ and $I(T)$, which comes from the fitted curve, with the corresponding cumulative number of infected individuals at times $\tau$ and $T$ in simulated data. We use Python's Scipy library to find optimal parameters for the logistic curves. Table \ref{s_shape_errors} reports the percent difference between the simulated and fitted values for different settings. We sum up the cumulative number of infected individuals in every district to obtain the values for Sierra Leone. Overall, we obtain good fits for the values of cumulative number of infected individuals. 
 
\begin{table}[H]
	\centering
	\caption{Percent difference between the simulated and fitted cumulative number of infected individuals at $\tau = 60$ and $T=188$.}
	\vspace{0.2cm}
	\large
	{\begin{tabular}{|c|c|c|c|}
			\hline
			Cumulative Infected & Setting 1 & Setting 2 & Setting 3 \\ \hline
			$I(\tau)$ & $22.56\%$ & $-3.98\%$ & $5.49\%$ \\ \hline
			$I(T)$ & $2.78\%$ & $13.72\%$ & $-0.37\%$ \\ \hline
	\end{tabular}}
	\label{s_shape_errors}
\end{table}

Recall that we estimate coefficients $a_n^1$, $b_n^1$, $a_n^2$ and $b_n^2$ in Equation \eqref{logistic_curves} for each population first for each stage. Then, in order to include the effect of initial infected population and resources allocated, we construct a regression for $K_n^1$ and $K_n^2$. To estimate the coefficients of these regressions, we generate data points by generating random initial states and resources. Specifically, the initial infected number of individuals are generated from the initial infected number of individuals observed from real data plus a normally distributed noise with mean zero and standard deviation 50 to obtain a wide range of values. Then, we find the initial number of susceptible individuals by subtracting the number of infected individuals from the total population in every district. We also generate random first and second stage capacity vectors for every initial state by using the real capacity vectors. We add a normally distributed random noise term to the total capacity with mean zero and standard deviation 50. Then we order the districts randomly and allocate the total number of beds by a uniform distribution between zero and the remaining capacity. The remaining capacity is found by subtracting random number of beds from the total available beds until all the districts are covered. The regression coefficients for all districts are found using lasso regression from the Scikit-learn library in Python. The average $R^2$ values for the two regression models are given in Table \ref{r_squared}. Our results yield high quality approximations for the cumulative number of infected individuals which we then use in the TSSP model. We provide the detailed analysis for the curve fit and regression model results in Appendix \ref{data_analysis}.

\begin{table}[H]
	\centering
	\caption{$R^2$ values for the regression models for Stages 1 and 2 for all settings.}
	\vspace{0.2cm}
	\large
	{\begin{tabular}{|c|c|c|c|}
			\hline
			& Setting 1 & Setting 2 & Setting 3 \\ \hline
			Stage 1 & 0.93 & 0.91 & 0.92 \\ \hline
			Stage 2 & 0.88 & 0.87 & 0.89 \\ \hline
	\end{tabular}}
	\label{r_squared}
\end{table}

\subsection{Comparison of Policies}
The goal of this section is to compare the results of policies produced by our optimization framework with that of what actually happened during the 2014 Ebola outbreak in Sierra Leone. The results can be interpreted as improvement that could have been achieved by implementing the proposed framework. Of course, these interpretations are based on the assumptions considered in this study. To that end, how we estimate the parameters of the model from available data from 2014 Ebola outbreak in Sierra Leone is explained next.

There are 1544 beds in ETUs across Sierra Leone \citep{EbolaBeds}. We assume that the 634 beds with an $operating$ status are available in the first stage, and the remaining 910 beds with a $pending$ status become available in the second stage. To estimate the opening costs we use data from \cite{Johnson2016} for the Western Area, and calculate the total construction cost per individual. Then, we use this information to calculate the ETU opening costs for the other districts proportional to their total population. The operating costs include medical and transportation costs. The medical costs are estimated by \cite{Bartsch2015}. They include treatment and personal protective equipment costs per patient in Sierra Leone for the 2014 Ebola epidemic. The essential items needed in ETUs are delivered from depots in Bombali, Kenema, and Port Loko, which constutite the transportation costs. We assume that every district is served by its closest depot, and then we calculate distance based costs considering the fuel prices. We used the fuel prices from Sierra Leone in 2015 \citep{fuel_prices}. Using this cost structure, the cost of bed allocations in Sierra Leone correspond to $\$505,000$ in the first stage and $\$470,000$ in the second stage. We consider three budgetary scenarios for the second stage: low, medium, and high. The $\$470K$ available for the second stage is used as the ``medium" value. The ``low" and ``high" values are set to $\$235K$ ($50\%$ less) and $\$705K$ ($50\%$ more), respectively. The probabilities of having low, medium, and high budgets are set to $0.25$, $0.5$, and $0.25$, respectively. The integrated TSSP model is solved using Gurobi 8.0. 

Table \ref{compare_results} reports the improvement due to optimized decisions compared to real data in terms of the cumulative number of infected individuals in each of the three settings. The following observations are in order. First, even with 20 days of data, our approach can lead to improvements of about $1.5\%$, which translates to averting $220$ Ebola cases. Second, if the decision maker uses the updated data when the extra budget is revealed and re-optimizes bed allocation decisions for the second stage, then the improvement can be almost doubled (i.e., 206 more cases can be averted). This suggests that value of information in a data-driven approach for infectious disease control is significant. Finally, a retrospective analysis shows that up to $5\%$ improvement is possible by following different allocations, which is the focus of our next discussion. 
Note that the improvements are calculated for a period of 188 days and one may expect larger improvements over longer periods.
In a data-driven setting, one may interpret the $5\%$ improvement as the value of perfect information. A similar conclusion can be drawn when comparing the $5\%$ improvement in Setting 3 with $1.5\%$ improvement in Setting 1 where we can quantify the value of information on the cumulative number of infected individuals.

\begin{table}[H]
	\centering
	\caption{Percent improvement due to optimization.}
	\vspace{0.2cm}
	\large
	{\begin{tabular}{|c|c|c|}
			\hline
			Setting 1 & Setting 2 & Setting 3 \\ \hline
			$1.49\%$ & $2.89\%$ & $4.37\%$ \\ \hline
	\end{tabular}}
	\label{compare_results}
\end{table} 

\begin{figure}[H]
	\centering
	\includegraphics[width=0.9\textwidth]{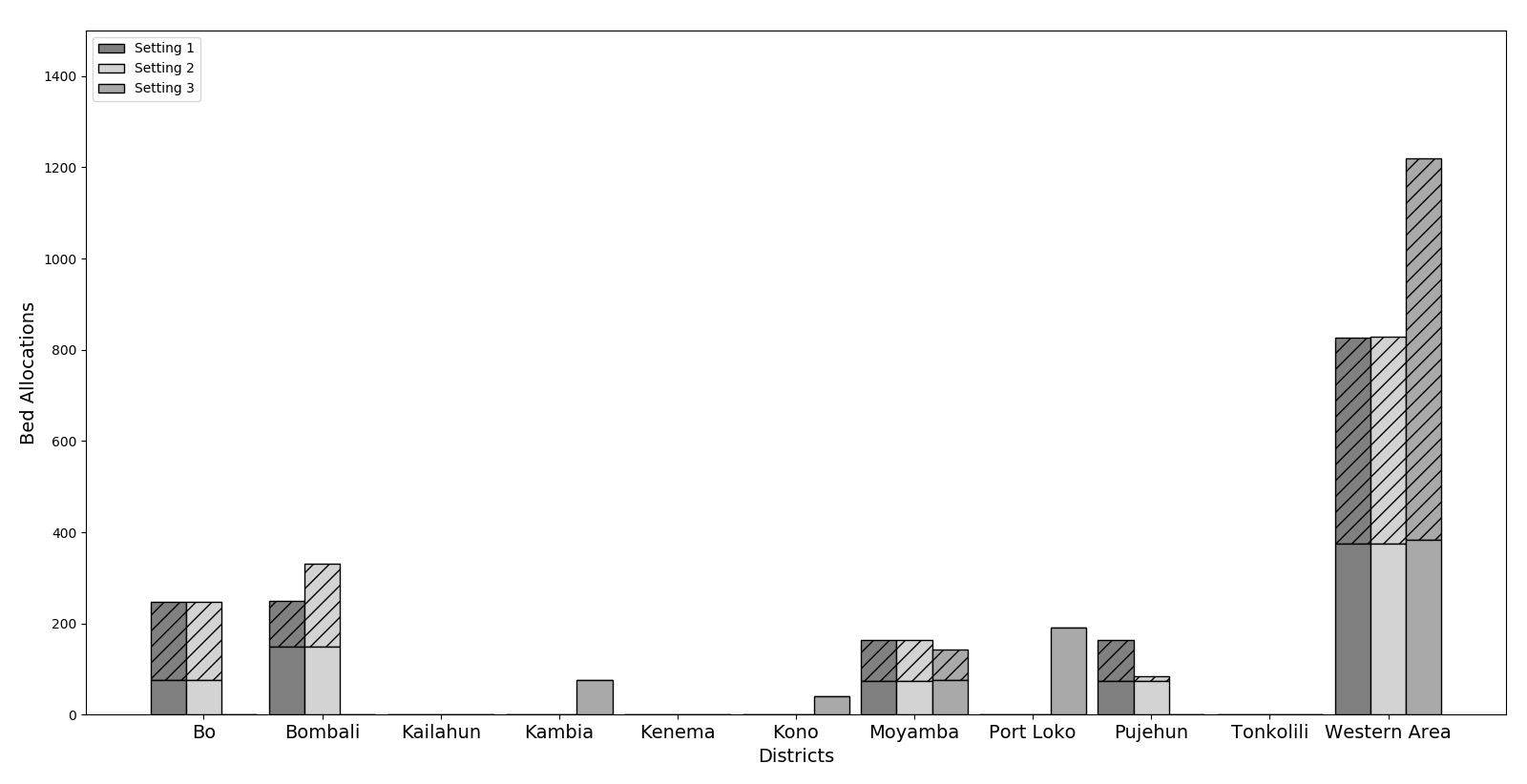}
	\caption{First and second stage bed allocations obtained by settings 1, 2, and 3.}
	\label{allocation_settings}
\end{figure}

The proposed optimization approach does lead to improvements in terms of reducing the number of infected individuals, but it is also interesting to analyze the decisions made by the proposed optimization and compare them to the actual decisions made during the 2014 Ebola epidemic. Figure \ref{allocation_settings} shows the bed allocations obtained by the optimization approach for the three settings. The lower and upper parts of the bars show the first and second stage bed allocations for every district, respectively. Recall that, the first stage decisions are the same for Settings 1 and 2, whereas the second stage decisions for Setting 2 are updated using new information. As can be seen in Figure \ref{allocation_settings}, under Setting 2 the resource allocation is shifted from Pujehun to Bombali. Initially, Bombali has low number of infected individuals. However, the epidemic grows fast in Bombali in the later time periods. Therefore, with the new information available, the resources are shifted to address this change in trends. In Setting 3, the resources are allocated in some districts that do not even have treatment units opened in Settings 1 and 2. In addition, in Setting 3, more resources are allocated to Western Area. As a result of this change in resource allocation, Setting 3 results in lower cumulative number of infected individuals in Sierra Leone. This observation suggests that the epidemic control strategies should prioritize more populated areas, while covering other districts with lower levels of resources.

\begin{figure}[H]
	\centering
	\includegraphics[width=\textwidth]{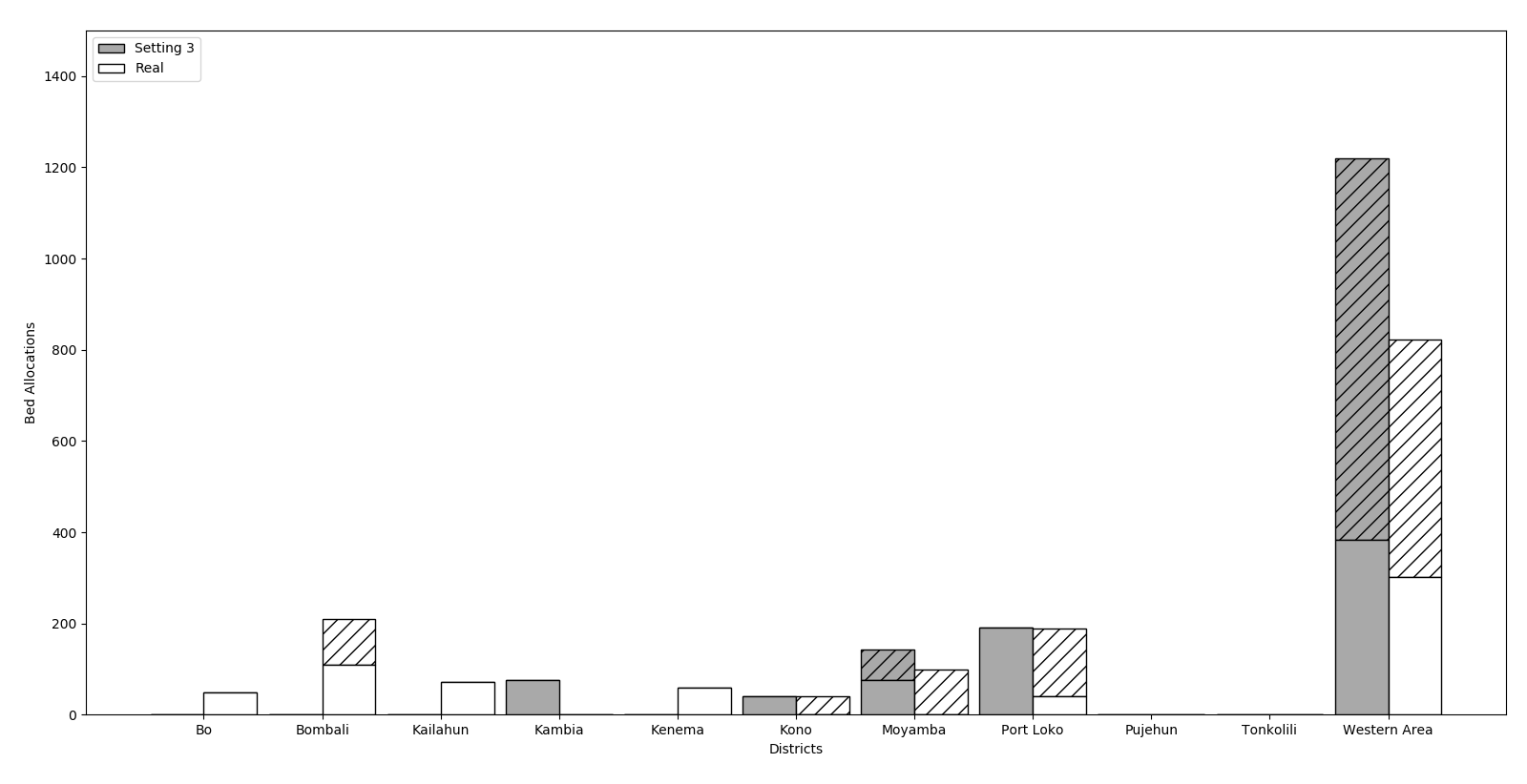}
	\caption{First and second stage bed allocations obtained by Setting 3 and real bed allocations.}
	\label{allocation_real}
\end{figure}

Figure \ref{allocation_real} compares the allocation strategies obtained by perfect information in Setting 3 to the real bed allocations. Similar to Figure \ref{allocation_settings}, the lower and upper parts of the bars correspond to first and second stage allocations, respectively. In both cases, Western Area receives the most allocation, but our optimization model allocates more to this district. Especially, in the second stage, when the growth of the epidemic in Western Area is faster than the other districts, our model takes appropriate precautions, and this results in a lower cumulative number of infected individuals by $11\%$ in Western Area compared to the cumulative number of infected individuals obtained by the real bed allocations. In Kono, Moyamba, and Port Loko, our model assigns the beds in the first stage as opposed to real bed allocations where more resources are allocated in the second stage. In Kono and Port Loko, our allocation results in a higher cumulative number of infected individuals. However, in Moyamba by allocating more resources, $128\%$ improvement in the cumulative number of infected individuals is achieved. In the real setting, the first stage budget is used in opening more treatment units in more districts. On the other hand, in our optimization setting, the budget is used in allocating more beds with fewer number of ETUs. 

\subsection{Sensitivity Analysis}\label{sensitivity}

We conduct a sensitivity analysis by changing the first stage budgets using Setting 3 with perfect information in our simulation model. In this analysis, we keep the set of second stage budget the same. The first stage budgets used are $\$250,000$, $\$375,000$, $505,000$, and $\$625,000$. Figure \ref{scenario_budget} shows the difference in the cumulative number of infected individuals with respect to different budgets.

\begin{figure}[H]
	\centering
	\begin{subfigure}{\textwidth}
		\centering
		\includegraphics[scale=0.85]{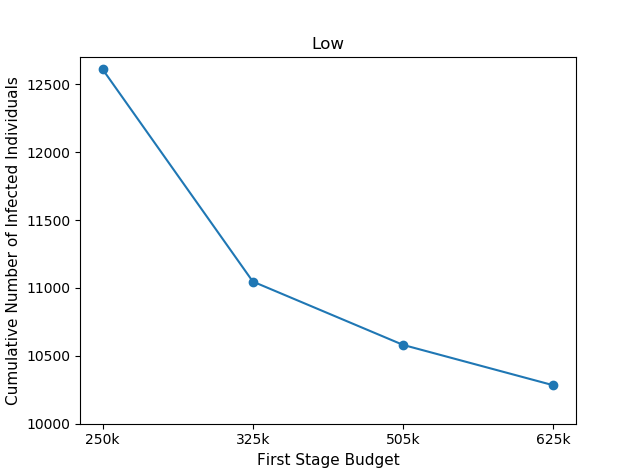}
		\caption{Low level extra budget.}
	\end{subfigure}
	\\
	\begin{subfigure}{\textwidth}
		\centering
		\includegraphics[scale=0.85]{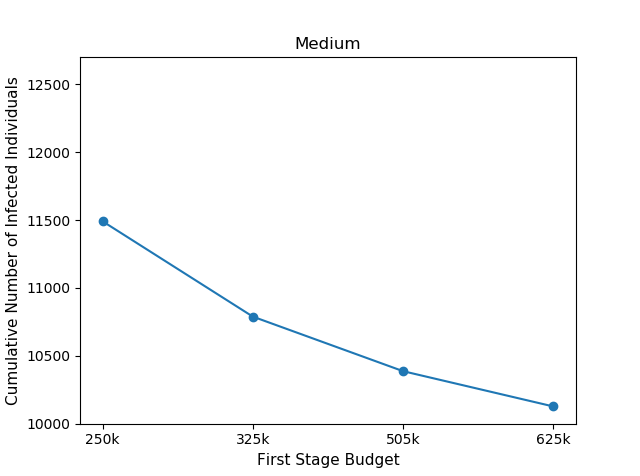}
		\caption{Medium level extra budget.}
	\end{subfigure}
\end{figure}

\begin{figure}[H]\ContinuedFloat
	\centering
	\begin{subfigure}{\textwidth}
		\centering
		\includegraphics[scale=0.85]{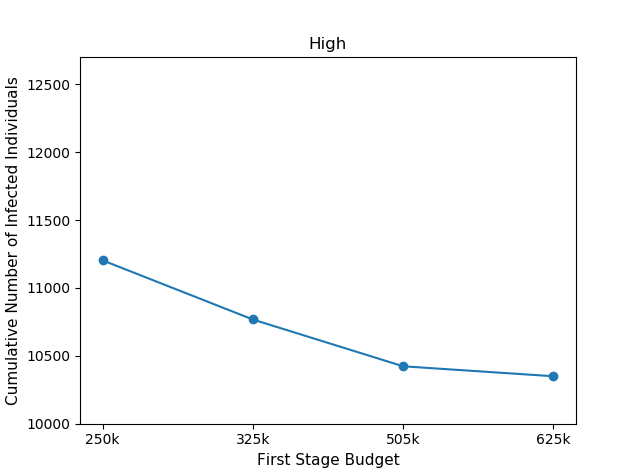}
		\caption{High level extra budget.}
	\end{subfigure}
	\caption{Cumulative number of infected individuals in Sierra Leone with respect to different first stage budgets.}
	\label{scenario_budget}
\end{figure}

As expected, when we increase the first stage budget, the number of cumulative infected individuals is decreased at time $T=188$ for every second stage extra budget level. However, the improvement is not linear. This is because after a critical level of resources is achieved, the infection progression slows down, and the additional resources have only marginal effect as the need for these resources has diminished. Therefore, depending on the first stage budgets, decision makers might choose to allocate some of these resources to other activities.

\vspace{-4mm}
\section{Conclusion}\label{conclusion}
In  this study we formulated the problem of allocating resources to control an infectious disease in a metapopulation where the budget is uncertain. We developed a TSSP model to decide opening of ETUs and bed allocations using an initial budget followed by a random budget becoming available in the second stage for additional bed allocations.

We included the disease dynamics in our TSSP formulation via a data-driven approach and using the observations that the number of cumulative infected individuals follows an \textit{s}-shaped curve, which resulted in a tractable model. The functional form of the cumulative number of infected individuals in our model depended on the initial number of infected individuals and the corresponding bed allocations, among other factors. We also observed, with extensive simulations, that it suffices to consider a linear regression model for only one of the parameters ($K_n$) of logistic curves in order to incorporate initial states and actions taken which resulted in a tractable model with high-quality solutions. 

To estimate the parameters, we developed a simulation model for the 2014 Ebola epidemic and adopted an online approach that mimics the real decision making process. When a decision maker is initially observing the epidemic progress, she has to use limited amount of data to estimate the future trajectories of the epidemic. Then, she solves an optimization problem to find ``optimal" actions. After receiving more observations about the epidemic and the resolution of the additional uncertain budget, the decision maker can update the optimization model parameters to find better actions for epidemic control. Our results showed that, if the decision makers used the proposed approach, there could have been 426 (i.e., about a 3\% improvement) averted cases in Sierra Leone which, while seems small, is significant. This approach allows for a flexible decision making methodology by using simulation models validated by even limited amount of data. We also provided a what-if  analysis, where we used all the available data, similar to other data-driven optimization models in the literature.

Our results indicated that the ETU opening and bed allocation decisions were affected by the amount of information that was available to the decision makers. When we updated the information after the additional budget is revealed, the first stage decisions were fixed and could not be altered. It is worth noting that our analysis showed completely different allocation strategies when we used all the available data, which resulted in the highest number of averted cases. However, when there is an emerging epidemic, decision makers rarely have access to perfect information. Therefore, our analysis quantifies the extent to which updating information can improve decision making in containing epidemics. 
  
When we compared our allocation policies to the real bed allocations in Sierra Leone, we observed that the Western Area, which is the most populated region of Sierra Leone, received most of the allocations in both cases. However, our policy allocated more to this district. Our policy also showed difference in the number of districts that had ETUs. In the real case, most of the budget was used in the initial opening cost of the treatment units, whereas our policy opens fewer ETUs but allocates more beds. These results provided us with valuable insights on how to allocate resources in an unknown limited budget setting to contain an epidemic. Overall, we achieved 645 (i.e., 5\%) averted cases by using the policy generated by our optimization model, compared to the real bed allocations, when we used all the available data. 

The framework of our model is flexible to be applied for other types of infectious diseases with different disease dynamics and control policies. Our simulation model provided detailed characteristics for the Ebola epidemic, however, other compartments could be included for disease specific purposes. In our model, we assumed a centralized decision maker, which has access to all the budget to be allocated for infectious disease control. One natural extension to our problem can be comparing coordinated and uncoordinated decision making processes of the relief organizations by adopting similar online decision making processes with information asymmetry.

\vspace{-4mm}
\bibliographystyle{apa}
\bibliography{references}


\appendix
\section{Supplementary Material for the Validated Simulation Model} \label{districts_validated}
In this section, we present the results for our validated simulation for every district of Sierra Leone. Note that, we report the cumulative number of infected individuals for the national level in Sierra Leone in Section \ref{result}. Table \ref{validation_appendix} shows the calibration and validation errors for every district corresponding to our three different settings. Recall that for Settings 1 and 2 we use 20 and 60 days of data respectively, and for Setting 3 we use all available data. In Figure \ref{time_series_validation}, we show the time series generated by the validated simulation and real data for every district for Setting 3. 

\begin{table}[h!]
\centering
\caption{Average percentage difference between real and simulated number of infected individuals for every district in Sierra Leone.}
{\begin{tabular}{|c|c|c|c|c|c|c|}
\hline
 & \multicolumn{2}{c|}{Setting 1} & \multicolumn{2}{c|}{Setting 2} & \multicolumn{2}{c|}{Setting 3} \\ \hline
Districts & Calibration & Validation & Calibration & Validation & Calibration & Validation \\ \hline
Bo & $21.07\%$ & $6.88\%$ & $12.44\%$ & $-6.85\%$ & $9.61\%$ & $-1.06\%$ \\ \hline
Bombali & $42.06\%$ & $38.94\%$ & $53.25\%$ & $29.66\%$ & $24.39\%$ & $-6.19\%$ \\ \hline
Kailahun & $-0.42\%$ & $-11.66\%$ & $3.76\%$ & $1.17\%$ & $1.84\%$ & $1.32\%$ \\ \hline
Kambia & $-35.01\%$ & $-1.63\%$ & $9.73\%$ & $24.02\%$ & $14.43\%$ & $-6.02\%$ \\ \hline
Kenema & $8.28\%$ & $4.65\%$ & $3.94\%$ & $1.22\%$ & $3.66\%$ & $2.58\%$ \\ \hline
Kono & $-406.26\%$ & $-73.62\%$ & $-195.66\%$ & $-27.95\%$ & $-135.45\%$ & $-3.3\%$ \\ \hline
Moyamba & $-16.41\%$ & $-7.93\%$ & $8.87\%$ & $7.7\%$ & $2.99\%$ & $-8.1\%$ \\ \hline
Port Loko & $29.83\%$ & $28.97\%$ & $31.54\%$ & $6.51\%$ & $17.48\%$ & $-11.14\%$ \\ \hline
Pujehun & $17.68\%$ & $-13\%$ & $11.08\%$ & $-0.05\%$ & $14.21\%$ & $0.55\%$ \\ \hline
Tonkilili & $20.35\%$ & $-13.95\%$ & $20.34\%$ & $-14.61\%$ & $12.81\%$ & $-1.97\%$ \\ \hline
Western Area & $33.44\%$ & $27.67\%$ & $33.69\%$ & $14.54\%$ & $20.97\%$ & $-2.94\%$ \\ \hline
\end{tabular}}
\label{validation_appendix}
\end{table}

\begin{figure}
\centering
\begin{subfigure}{\textwidth}
\centering
\includegraphics[scale=0.6]{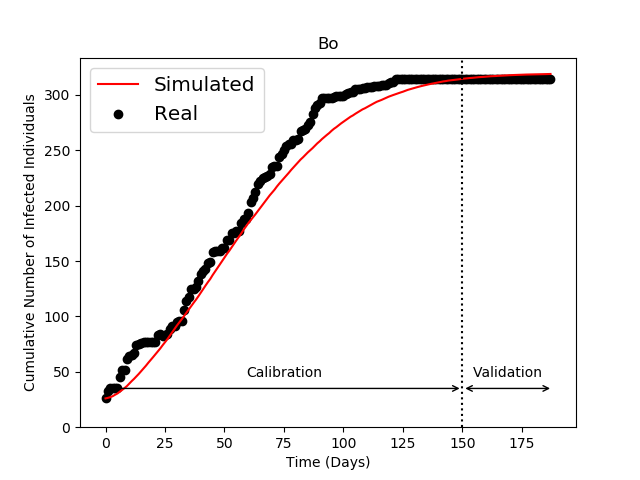}
\end{subfigure}
\\
\begin{subfigure}{\textwidth}
\centering
\includegraphics[scale=0.6]{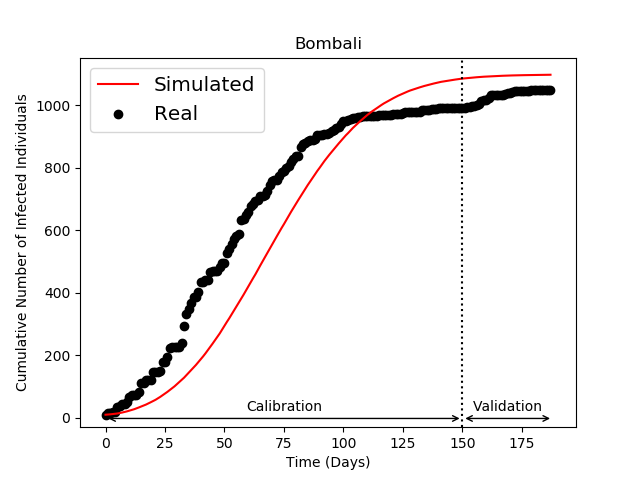}
\end{subfigure}
\\
\begin{subfigure}{\textwidth}
\centering
\includegraphics[scale=0.6]{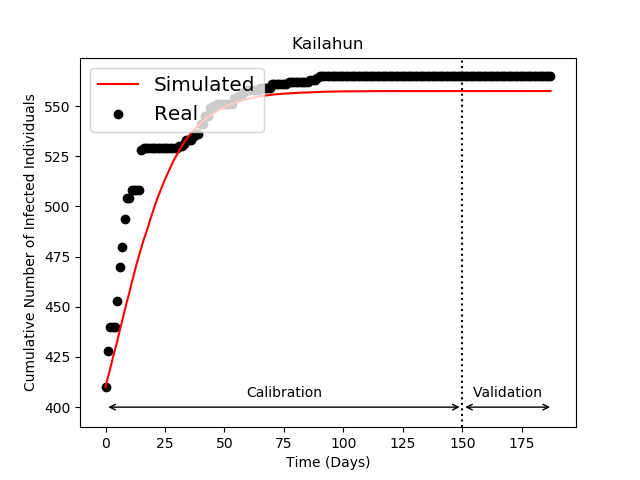}
\end{subfigure}
\end{figure}

\begin{figure}\ContinuedFloat
\centering
\begin{subfigure}{\textwidth}
\centering
\includegraphics[scale=0.6]{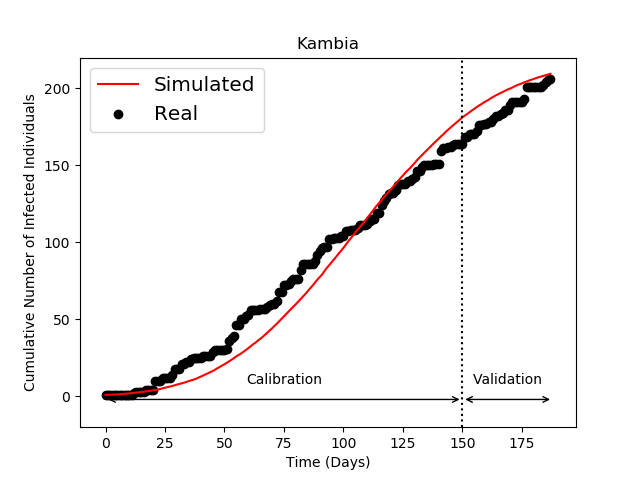}
\end{subfigure}
\\
\begin{subfigure}{\textwidth}
\centering
\includegraphics[scale=0.6]{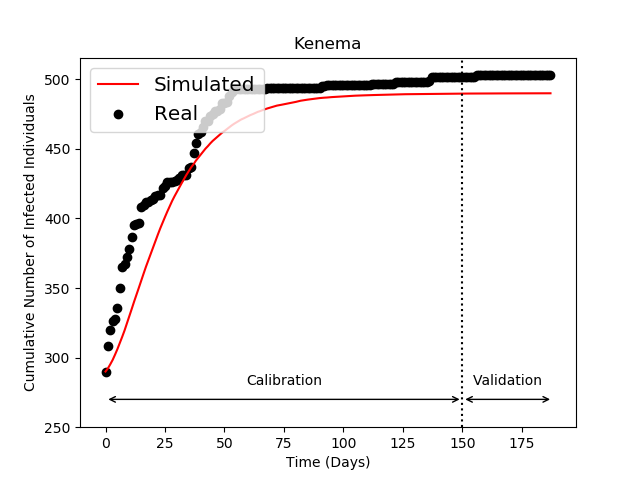}
\end{subfigure}
\\
\begin{subfigure}{\textwidth}
\centering
\includegraphics[scale=0.6]{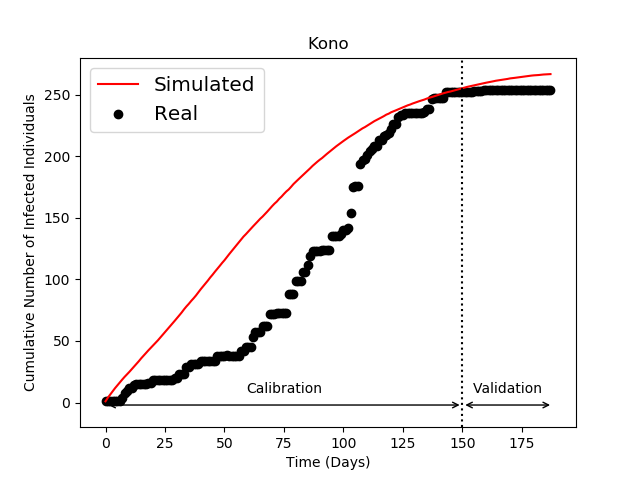}
\end{subfigure}
\end{figure}

\begin{figure}\ContinuedFloat
\centering
\begin{subfigure}{\textwidth}
\centering
\includegraphics[scale=0.6]{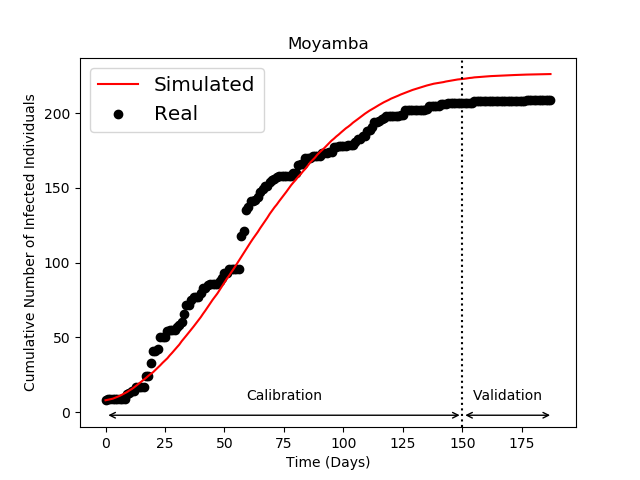}
\end{subfigure}
\\
\begin{subfigure}{\textwidth}
\centering
\includegraphics[scale=0.6]{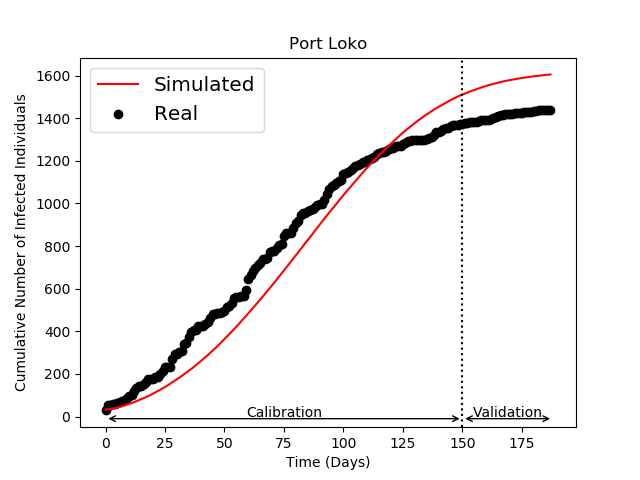}
\end{subfigure}
\\
\begin{subfigure}{\textwidth}
\centering
\includegraphics[scale=0.6]{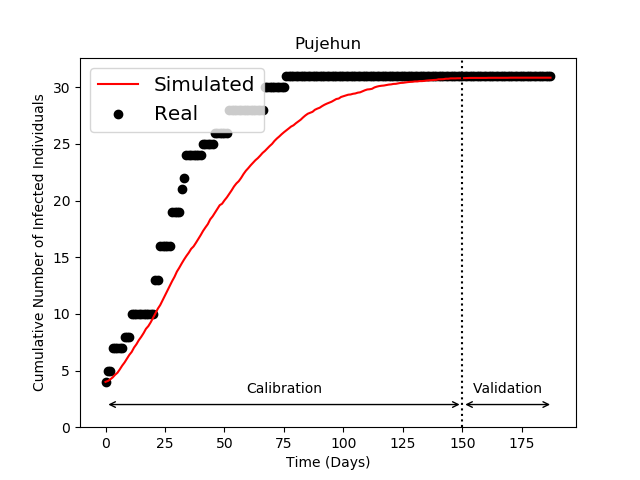}
\end{subfigure}
\end{figure}

\begin{figure}\ContinuedFloat
\centering
\begin{subfigure}{\textwidth}
\centering
\includegraphics[scale=0.6]{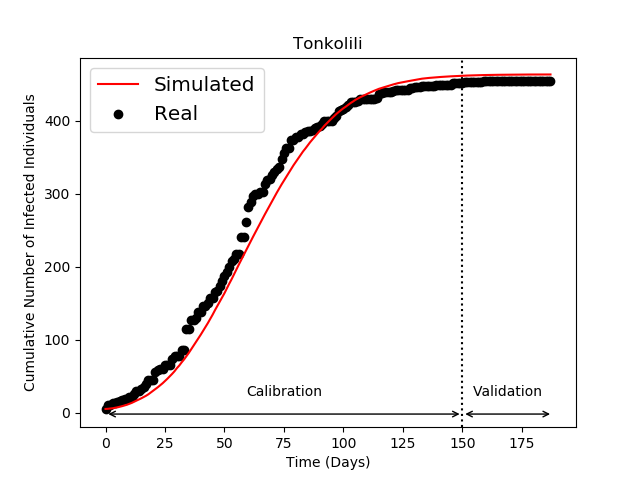}
\end{subfigure}
\\
\begin{subfigure}{\textwidth}
\centering
\includegraphics[scale=0.6]{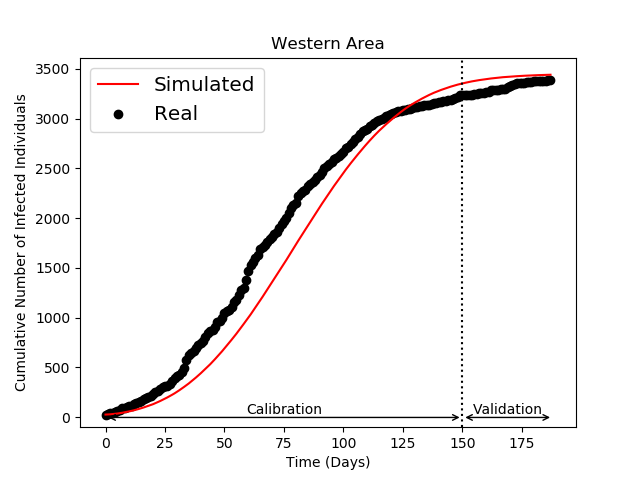}
\end{subfigure}
\caption{Cumulative number of infected individuals for every district in Sierra Leone showing the real and simulated data.}
\label{time_series_validation}
\end{figure}

\newpage
\section{Supplementary Material for Data Analysis} \label{data_analysis}
In this section, we provide the details of the data analysis carried out for estimating the parameters of the \textit{s}-shaped curves for every district as well as regression model assumptions. 

As our first task, we estimate the parameters of the \textit{s}-shaped curves of the form $\frac{K_n}{e^{-(a_n  t + b_n)}}$ that represent the cumulative number of infected individuals. We divide the data set into two; from time 0 to $\tau$, and from $\tau$ to $T$. We assume that $\tau = 60$, and $T=188$ in this study. For every district we generate the time series by our validated simulation model and fit an \textit{s}-shaped curve by minimizing the sum of squared errors. Note that we have two different curves for the data set that represents the two stages of our optimization model. Figure \ref{curves_districts} shows the fitted and real cumulative number of infected individuals in Setting 3, and Table \ref{difference_all} presents the percent difference of $I(\tau)$ and $I(T)$ for our three settings for every district of Sierra Leone.

\begin{figure}
\centering
\begin{subfigure}{\textwidth}
\centering
\includegraphics[scale=0.6]{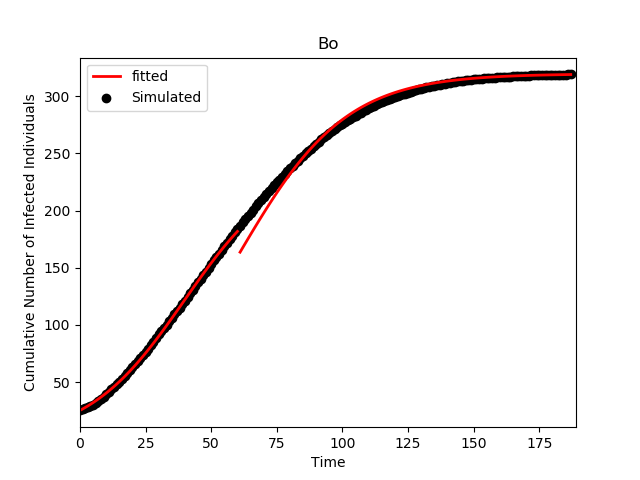}
\end{subfigure}
\\
\begin{subfigure}{\textwidth}
\centering
\includegraphics[scale=0.6]{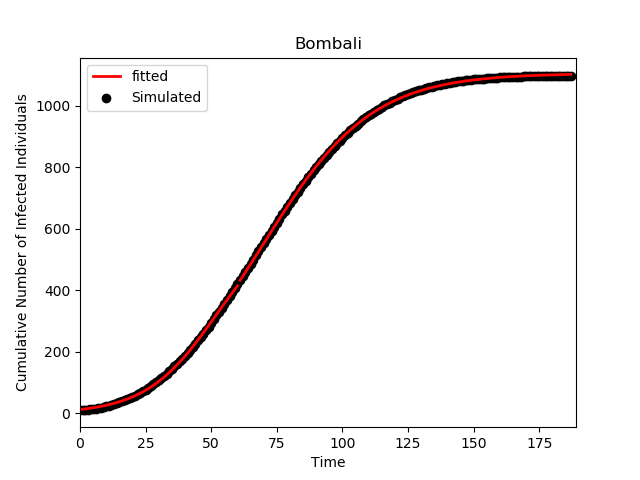}
\end{subfigure}
\\
\begin{subfigure}{\textwidth}
\centering
\includegraphics[scale=0.6]{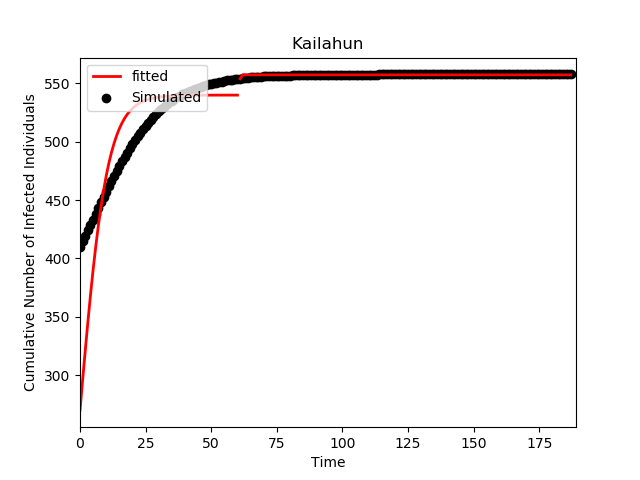}
\end{subfigure}
\end{figure}

\begin{figure}\ContinuedFloat
\centering
\begin{subfigure}{\textwidth}
\centering
\includegraphics[scale=0.6]{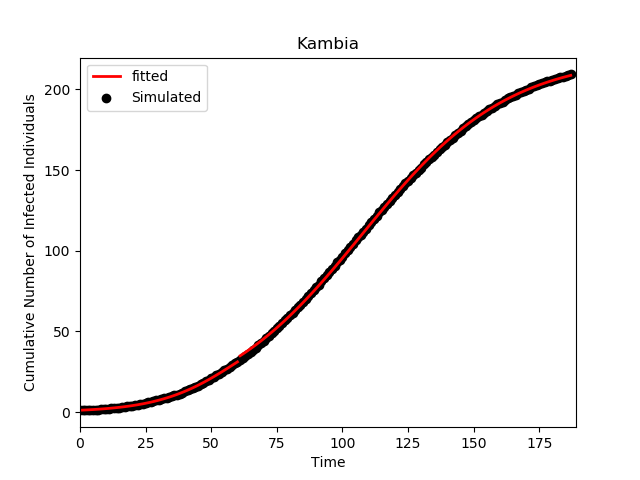}
\end{subfigure}
\\
\begin{subfigure}{\textwidth}
\centering
\includegraphics[scale=0.6]{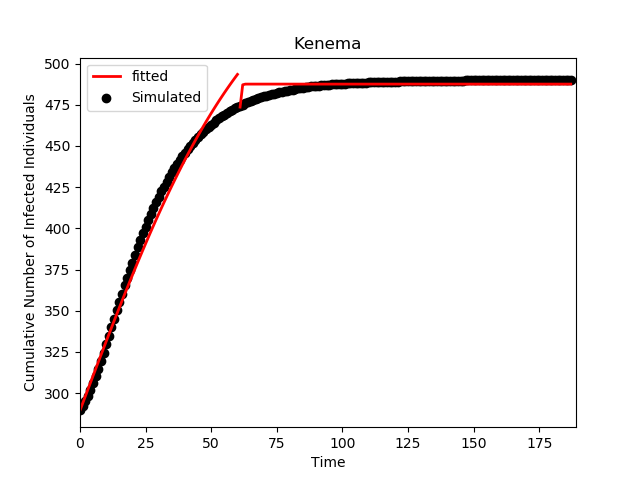}
\end{subfigure}
\\
\begin{subfigure}{\textwidth}
\centering
\includegraphics[scale=0.6]{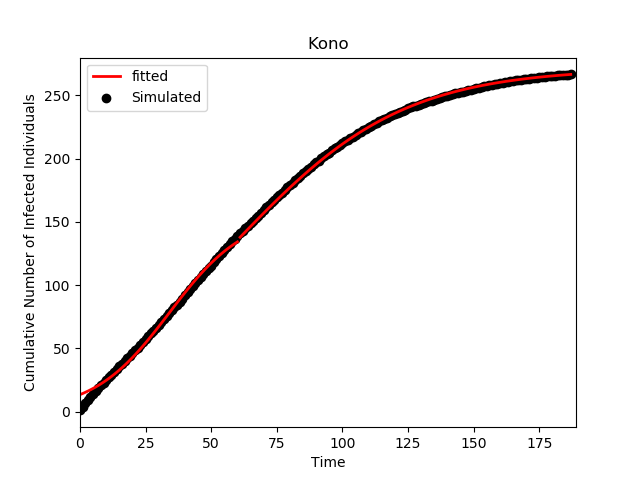}
\end{subfigure}
\end{figure}

\begin{figure}\ContinuedFloat
\centering
\begin{subfigure}{\textwidth}
\centering
\includegraphics[scale=0.6]{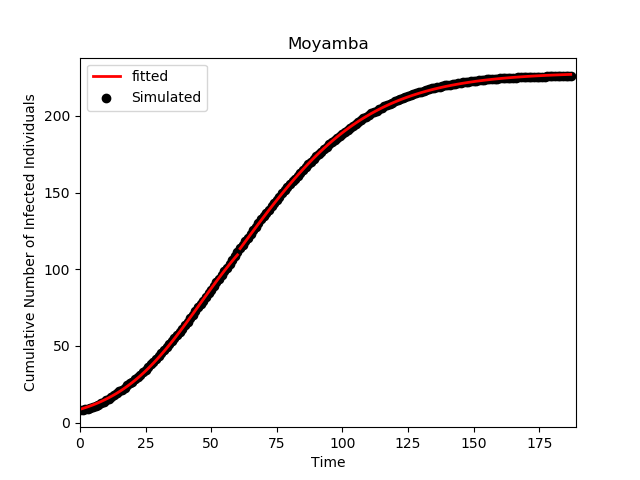}
\end{subfigure}
\\
\begin{subfigure}{\textwidth}
\centering
\includegraphics[scale=0.6]{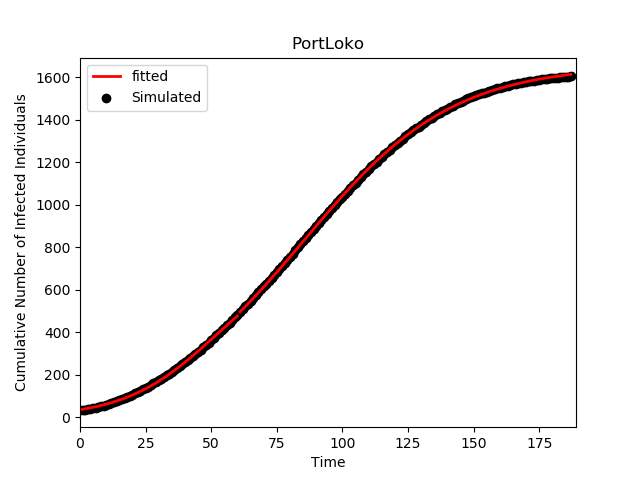}
\end{subfigure}
\\
\begin{subfigure}{\textwidth}
\centering
\includegraphics[scale=0.6]{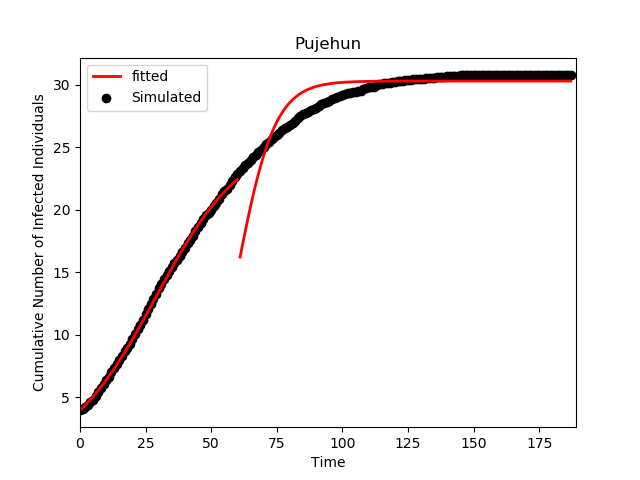}
\end{subfigure}
\end{figure}

\begin{figure}\ContinuedFloat
\centering
\begin{subfigure}{\textwidth}
\centering
\includegraphics[scale=0.6]{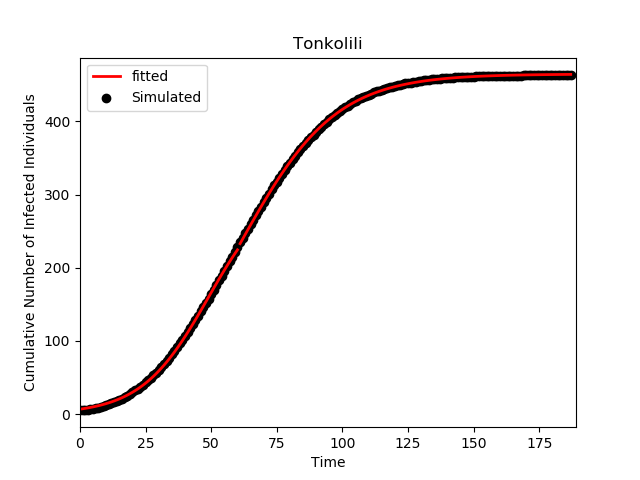}
\end{subfigure}
\\
\begin{subfigure}{\textwidth}
\centering
\includegraphics[scale=0.6]{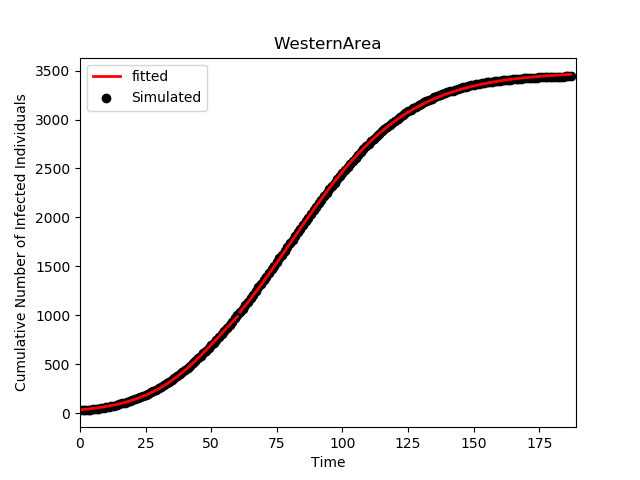}
\end{subfigure}
\caption{Cumulative number of infected individuals with respect to fitted logistic curves and simulated data.}
\label{curves_districts}
\end{figure}

\begin{table}[h!]
\centering
\caption{Percent difference for fitted and simulated data for $I(\tau)$ and $I(T)$ for the three settings for all districts of Sierra Leone.}
{\begin{tabular}{|c|c|c|c|c|c|c|}
\hline
 & \multicolumn{2}{c|}{Setting 1} & \multicolumn{2}{c|}{Setting 2} & \multicolumn{2}{c|}{Setting 3} \\ \hline
Districts & $I(\tau)$ & $I(T)$ & $I(\tau)$ & $I(T)$ & $I(\tau)$ & $I(T)$ \\ \hline
Bo & $2.92\%$ & $-0.24\%$ & $1.18\%$ & $-0.028\%$ & $1.32\%$ & $0.03\%$ \\ \hline
Bombali & $37.26\%$ & $1.56\%$ & $19.95\%$ & $0.27\%$ & $12.55\%$ & $-0.37\%$ \\ \hline
Kailahun & $-0.28\%$ & $0.66\%$ & $2.46\%$ & $0.03\%$ & $2.5\%$ & $0.03\%$ \\ \hline
Kambia & $21.56\%$ & $1.04\%$ & $21.09\%$ & $0.83\%$ & $15.41\%$ & $0.28\%$ \\ \hline
Kenema & $-0.25\%$ & $4.42\%$ & $-1.84\%$ & $0.51\%$ & $-4.2\%$ & $0.48\%$ \\ \hline
Kono & $2.02\%$ & $-0.02\%$ & $1.84\%$ & $-0.01\%$ & $3.32\%$ & $0.01\%$ \\ \hline
Moyamba & $7.29\%$ & $-0.32\%$ & $1.76\%$ & $-0.38\%$ & $1.67\%$ & $-0.46\%$ \\ \hline
Port Loko & $17.19\%$ & $0.77\%$ & $7.16\%$ & $-0.04\%$ & $1.51\%$ & $-0.56\%$ \\ \hline
Pujehun & $1.66\%$ & $-0.04\%$ & $1.09\%$ & $0.18\%$ & $1.81\%$ & $1.7\%$ \\ \hline
Tonkolili & $3.68\%$ & $1.31\%$ & $2.91\%$ & $1.24\%$ & $10.88\%$ & $-0.26\%$ \\ \hline
Western Area & $34.05\%$ & $3.81\%$ & $15.8\%$ & $-0.47\%$ & $12.94\%$ & $-0.58\%$ \\ \hline
\end{tabular}}
\label{difference_all}
\end{table}

After finding the initial parameters of the logistic curves, we fit a regression model for $K_n^1$ and $K_n^2$ and fix the other parameters. Our observations show that these parameters are functions of the bed capacities and the initial number of infected individuals. Therefore, we randomly generate initial states and bed capacities for the first and second stages, then simulate the system to obtain the time series for cumulative number of infected individuals with random initial settings. After that, we fit a logistic curve to these simulated time series to obtain data to be used in our regression model. As the number of independent variables is large, we use a lasso regression model. Table \ref{r_squared_districts} show the $R^2$ values of the two regression models for the two logistic curve parameter $K_n$ for all districts of Sierra Leone. 

\begin{table}[h!]
\centering
\caption{$R^2$ values for of regression models for each stage and district.}
{\begin{tabular}{|c|c|c|c|c|c|c|}
\hline
Training Time & \multicolumn{2}{c|}{$t=20$} & \multicolumn{2}{c|}{$t=60$} & \multicolumn{2}{c|}{$t=150$} \\ \hline
District & Stage 1 & Stage 2 & Stage 1 & Stage 2 & Stage 1 & Stage 2 \\ \hline
Bo & 0.97 & 0.97 & 0.96 & 0.93 & 0.97 & 0.89 \\ \hline
Bombali & 0.99 & 0.82 & 0.99 & 0.98 & 0.99 & 0.97 \\ \hline
Kailahun & 0.96 & 0.85 & 0.93 & 0.74 & 0.93 & 0.76 \\ \hline
Kambia & 0.95 & 0.98 & 0.75 & 0.99 & 0.76 & 0.98 \\ \hline
Kenema & 0.94 & 0.81 & 0.94 & 0.65 & 0.94 & 0.62 \\ \hline
Kono & 0.80 & 0.38 & 0.89 & 0.57 & 0.96 & 0.85 \\ \hline
Moyamba & 0.93 & 0.98 & 0.93 & 0.98 & 0.93 & 0.99 \\ \hline
Port Loko & 0.99 & 0.98 & 0.99 & 0.98 & 0.99 & 0.96 \\ \hline
Pujehun & 0.82 & 0.98 & 0.97 & 0.78 & 0.96 & 0.81 \\ \hline
Tonkolili & 0.95 & 0.98 & 0.76 & 0.99 & 0.77 & 0.99 \\ \hline
Western Area & 0.95 & 0.90 & 0.95 & 0.99 & 0.94 & 0.99 \\ \hline
\end{tabular}}
\label{r_squared_districts}
\end{table}
\listoffigures

\end{document}